\documentclass[english, 11pt]{article}

\usepackage[english]{babel}
\usepackage[latin1]{inputenc}
\usepackage{amsmath,amssymb,enumerate}
\usepackage{ulem}
\usepackage{setspace}
\usepackage{psfig,labelfig}
\usepackage{graphicx}
\usepackage[dvips]{epsfig}
\usepackage[T1]{fontenc}

\newtheorem{thm}{Theorem}[section]
\newtheorem{defi}[thm]{Definition}

\newtheorem{lem}[thm]{Lemma}
\newtheorem{con}[thm]{Conjecture}

\newcommand{\C}{{\mathbb C}}

\newcommand{\F}{{\mathbb F}}

\newcommand{\R}{{\mathbb R}}

\newcommand{\Z}{{\mathbb Z}}
\newcommand{\N}{{\mathbb N}}
\newcommand{\Zc}{{\mathcal Z}_g}
\newcommand{\dint} {\displaystyle\int}

\DeclareMathOperator{\vol}{vol}
\DeclareMathOperator{\sta}{Fix}

\title{A new lower bound for Hermite's constant for symplectic lattices}

\author{Muetzel, Bjoern\thanks{E-mail address : bjorn.mutzel@gmail.com}\\
 \\
\small Department of Mathematics, Universit\'e Montpellier 2, Case courrier 051, \\
\small 34095 Montpellier cedex 5, France \\[-0.8ex]}

\begin{document}

\maketitle

\begin{abstract}
In section 1 we give an improved lower bound on Hermite's constant $\delta_{2g}$ for symplectic lattices in even dimensions ($g=2n$) by applying a mean-value argument from the geometry of numbers to a subset of symmetric lattices. Here we obtain only a slight improvement. However, we believe that the method applied has further potential. In section 2 we present new families of highly symmetric (symplectic) lattices, which occur in dimensions of powers of two. Here the lattices in dimension $2^n$ are constructed with the help of a multiplicative matrix group isomorphic to $({\Z_2}^n,+)$. We furthermore show the connection of these lattices with the circulant matrices and the Barnes-Wall lattices.
\end{abstract}

\section{A new lower bound for Hermite's constant for symplectic lattices}

A \textit{lattice} $\Lambda$ of dimension $g$ is a discrete subgroup of $\R^g$ that spans $\R^g$. If $A$ is a matrix representation of a basis of $\Lambda$, then $\det(\Lambda)$, the \textit{determinant of $\Lambda$} is equal to $|\det(A)|$, the absolute value of the determinant of the matrix $A$. \textit{Hermite's constant or invariant} $\gamma_g$ is given by
\[
    \gamma_g = \max \limits_{\det(\Lambda)=1} \min \limits_{\lambda \in \Lambda \backslash \{0\}} \| \lambda \|_2 ^2.
\]
It is the maximal value that the squared norm of the shortest non-zero lattice vector - or systole - can attain among all lattices of determinant $1$.
We call a real $g$ dimensional lattice $\Lambda$ \textit{symmetric}, if there is an $O \in O_g(\R) \backslash \{\pm Id_g\}$, such that
\[
     O \Lambda = \Lambda.
\]
If $\Lambda=A\cdot \mathbb{Z}^g$, where $A$ is a matrix, whose columns form a basis of $\Lambda$, then this means that the matrix $O$ induces a change of basis. Therefore the above equality is equal to
\begin{equation}
O\cdot A = A\cdot R,
\text{ \ \ where } R \in GL_g(\mathbb{Z}).
\label{eqn:sym}
\end{equation}
The motivation to investigate symmetric lattices is the following result by Vorono\"i: \\
Among all lattices of given determinant, a lattice is called \textit{extreme}, if the length of its shortest non-zero lattice vector is a local maximum. Vorono\"i showed in \cite{vo} that this definition is equivalent to the lattice being \textit{perfect} and \textit{eutactic}. Especially, for a lattice to be extreme it has to have many systoles and often these lattices exhibit a high degree of symmetry.\\
The question that arises here is, whether symmetric lattices tend to have a larger systole in general than the average lattice and if this can be used to improve the lower bound for Hermite's constant.\\
\\
We pursue this idea in case of symplectic lattices. We give an improved lower bound for Hermite's constant for symplectic lattices $\delta_{2g}$ in even dimensions ($g=2n$) by applying a mean-value argument from the geometry of numbers to a subset of symmetric symplectic lattices. This approach applies also to real lattices in general. However in that case better bounds have been established in \cite{ro}.\\
\\
Here a \textit{symplectic (real) lattice} $\Lambda$ of dimension $2g$ is a lattice that has a basis with matrix representation in $Sp_{2g}(\R)$. $Sp_{2g}(\R)$ is  the \textit{symplectic group} over $\R$.
\[
Sp_{2g}(\mathbb{R})=\{R\in M_{2g}(\mathbb{R}) \mid \text{ }R^{T}\left(
\begin{array}{cc}
0 & Id_g \\
-Id_g & 0
\end{array}
\right) R=\left(
\begin{array}{cc}
0 & Id_g \\
-Id_g & 0
\end{array}
\right) \}.
\]
Let
\[
\mathcal{H}_g=\left\{  Z=X+iY \in M_g(\C) \mid X=X^T, Y=Y^T \text{ \ and  \ }Y > 0 \right\}
\]
be \textit{Siegel's upper half space}.\\
It was shown in \cite{bs} that a $2g$-dimensional lattice $\Lambda$ is symplectic if and only if there is a $Z \in \mathcal{H}_g$, such that $ \Lambda = P_Z \cdot \Z^{2g}$, where
\begin{equation}
P_{Z}=P_{X,Y}=\left(
\begin{array}{cc}
\sqrt{Y^{-1}} & 0 \\
0 & \sqrt{Y}
\end{array}
\right) \left(
\begin{array}{cc}
Id_g & X \\
0 & Id_g
\end{array}
\right)=
\left(
\begin{array}{cc}
\sqrt{Y^{-1}} & \sqrt{Y^{-1}}X \\
0 & \sqrt{Y}
\end{array}
\right).
\label{eq:PZ}
\end{equation}
In this way the symplectic lattices can be identified with the principally polarized abelian varieties (see \cite{bs}) and $\mathcal{H}_g$ is a \textit{moduli} or \textit{parameter space} for symplectic lattices. \\

In this setting Hermite's invariant $\delta_{2g}$ for symplectic lattices of dimension $2g$, is given by
\[
\delta_{2g}= \max\limits_{Z \in \mathcal{H}_g} \min\limits_{l \in \Z^{2g} \backslash \{0\}} \| P_Z \cdot l \|_{2}^2.
\]
It is the maximal value for Hermite's constant among the symplectic lattices. Note that by (\ref{eq:PZ}) the determinant of a symplectic lattice is always equal to $1$.\\
Explicit examples of symplectic lattices, where the lengths of the systoles of the corresponding lattices are a local maxima are given in \cite{ber1},\cite{ber} and \cite{bs}. Buser and Sarnak show in \cite{bs}:
\[
\frac{1}{\pi }\sqrt[g]{2 \cdot g!} \leq \delta_{2g}.
\]

We now give a slightly improved bound in the case, where the dimension is a multiple of four using a subset of symmetric lattices.
\begin{thm}
If $g=2n$, then
\[
 \frac{1}{\pi }\sqrt[g]{4 \cdot g!} \leq \delta_{2g}.
\]
\label{thm:Hermite_mu}
\end{thm}
The proof of the theorem is similar to the proof given in \cite{bs}. We present families of highly symmetric lattices in section 2. We furthermore show the connection of these lattices with the Barnes-Wall lattices (see \cite{bw} and \cite{ni}) and the circulant matrices (see \cite{da}). In the first case, we show that there exist matrix versions of the Barnes-Wall lattices that share certain symmetries with the matrices from which our families are derived. These matrices belong to the class of generalized circulant matrices or patterned matrices presented in \cite{ch} or \cite{di}. The families of lattices presented in section 2 are related to the symmetric lattices from the proof below. \\
\\
\textbf{proof} Let $Z=X+iY \in \mathcal{H}_{g}$. For a positive real number $r$ and a symplectic lattice  $\Lambda = P_Z\cdot \Z^{2g}$ we define $n_{r}(\Lambda)$ as the number of non-zero lattice vectors of $\Lambda$ whose lengths are smaller than or equal to $r$:

\[
n_{r}(\Lambda)=\#\{\lambda \in \Lambda \backslash \{0\}\text{ }|\text{ } \|\lambda\|_2  \leq r\}\text{.}
\]

For a $Z \in \mathcal{H}_g$ we obtain
\[
n_{r}(P_Z \cdot  \Z^{2g})=\sum_{l\in \mathbb{Z}^{2g}\backslash \{0\}}\chi
_{r^{2}}(l^{T}Q_{Z}l)\text{,}
\]

where $Q_Z = {P_{Z}}^T P_{Z}$ and $\chi _{r^{2}}$ is the characteristic function of the interval $[0,r^{2}]$.
We now consider the matrices $Z \in \mathcal{H}_g =\mathcal{H}_{2n}$ of the form
\[
Z=X+i\frac{1}{y^2}Id_{g}\text{,}
\]
where $y>0$, $X=X^{T}$, $KXK=X$ and $K$ is the orthogonal matrix

\[
K=\left(
\begin{array}{cccc}
0 & ... & 0 & 1 \\
... &  & -1 & 0 \\
0 & 1 & ... &   \\
-1 & 0 & ... & 0 \\
\end{array}
\right) \in M_{g}(\mathbb{R)}\text{.}
\]
The matrix $KXK=X={(X_{i,j})}_{i,j=1,..,g}$ has the following shape:

\begin{eqnarray*}
     X_{i,j}&=& X_{j,i} \\
     X_{i,j}&=& \left\{ {\begin{array}{*{20}c}
             {-X_{g+1-i,g+1-j}  \text{ \ if \ }    i+j \text{ is even }}\\
             {X_{g+1-i,g+1-j}  \text{ \ if \ }    i+j \text{ is odd }}\\
              \end{array} } \right. \text{.}
\end{eqnarray*}
This is equivalent to saying that the matrix $\left| X \right|$ is symmetric with respect to the first and second main diagonal and the diagonals below the second main diagonal of $X$ are alternatively the positive and negative image of the corresponding diagonals above the second main diagonal with respect to the reflection along this diagonal. The matrix $X$ is completely determined by the $\frac{g}{4}(g+2)$ entries $\{ X_{i,j} \mid (i \le j) \wedge (i \le g+1-j)\}$.\\
Under these conditions the matrix $P_Z=P_{X,\frac{1}{y^2}Id_g}$ has the following property:

\[
K'\cdot P_Z=P_Z\cdot K', \text{ \ \ where \ \ } K'=\left( {\begin{array}{*{20}c}
   K & 0  \\
   0 & K^{T}  \\
 \end{array} } \right).
\]

It follows from equation (\ref{eqn:sym}) that the corresponding lattice $P_Z \cdot \mathbb{Z}^{2g}$ is closed under the symmetry induced by the orthogonal matrix $K'$, whose only fixed point is zero. It follows that the lattice vectors $P_Z \cdot l$ and $P_Z \cdot K'l$ of the lattice $P_Z \cdot \mathbb{Z}^{2g}$ have the same length and $K'l=l$ if and only if $l=0$. Therefore there exist always 4 disjoint non-zero lattice vectors in $P_Z \cdot \mathbb{Z}^{2g}$ of the same length.\\

We now consider the set $V_g$ for $g=2n$:
\begin{eqnarray}
V_g &=&\{X\in M_{g}(\mathbb{R)}\text{ }\mid\text{ }X=X^{T},KXK=X,X_{ij}\in [0,1]\}\text{ \ \
and} \nonumber \\
\dim(V_g) &=&v_g=\frac{g(g+2)}{4}\text{.} 
\label{eq:Vg}
\end{eqnarray}
Here $\vol(V_g)=1$ as we integrate over the interval $[0,1]$ in each variable.
For a fixed $y$ we obtain the mean value $I(y)$ of $n_{r^{2}}(A_{Z},H_{Z})$ over the compact set $V_g$:
\[
I(y)=\frac{1}{\vol(V_g)}\dint\limits_{V_g}\text{ \ }\sum\limits_{l\in \mathbb{Z}%
^{2g}\backslash \{0\}}f(P_{X,\frac{1}{y^2}Id_{g}}\cdot l)\text{ \ }%
d X \text{.}
\]
Here $f:\R^{2g}\rightarrow \R$ is the function\ $f(x)=\chi_{r^{2}}(x^{T}\cdot x)$ and $d X=\prod\limits_{(i\leq j) \wedge (i \le g+1-j)}d X_{ij}$. With the help of the following \textbf{Lemma \ref{thm:mean_value}} and the volume formula for unit balls in dimension $2g$ we obtain

\begin{equation*}
\lim\limits_{y\rightarrow +\infty }I(y) =\int\limits_{\R^{2g}} { f(x) \,dx }  =\sigma _{2g} \cdot r^{2g} =\frac{\pi ^{g}\cdot r^{2g}}{g!}\text{.}
\end{equation*}

For $r^{2}<\frac{1}{\pi }\sqrt[g]{4 \cdot g!}=c$ one has

\begin{equation*}
\lim\limits_{y\rightarrow +\infty }I(y)<4.
\end{equation*}
As $I(y)$ is the mean number of non-zero lattice vectors with squared length smaller than or equal to $r^{2}$ in
$V_g$ and as $I(y)$ is decreasing for increasing values of $y$, there exists for a $r^{2}<c$ a $y^{\prime }$, such that the mean value $I(y^{\prime })$ is smaller than $4$. As in our case the number of non-zero lattice
vectors of fixed length is always a multiple of four, there exists an $X^{\prime }$ in $V_g$, such that the PPAV with period matrix $(X^{\prime },\frac{1}{{y^{\prime }}^2}Id_{g})$ for $r^{2}<c$ has no non-zero lattice vector with squared length smaller than $r^{2}$. Therefore we have found a lower bound for $\delta_{2g}$ and we obtain 
\[
\delta_{2g} \geq \frac{1}{\pi }\sqrt[g]{4 \cdot g!}. \ \ \ \square
\]

\begin{lem}
Let $g=2n$ and $f:\mathbb{R}^{2g}\rightarrow \mathbb{R}$ an integrable function
of compact support. Consider the function $I_{f}:\mathbb{R}%
^{+}\rightarrow \mathbb{R}$ which is defined by
\[
I_{f}(y)=\dint\limits_{V_g}\text{ \ }\sum\limits_{l\in \mathbb{Z}%
^{2g}\backslash \{0\}}f(P_{X,\frac{1}{y^2}Id_{g}}\cdot l)\text{ \ }%
d X\text{,}
\]
such that $P_Z$ and $V_g$ as in equation (\ref{eq:PZ}) and (\ref{eq:Vg}) respectively and $d
X=\prod\limits_{(i\leq j),\\ (i\leq g+1-j) }d X_{ij}$. One has:
\[
\lim\limits_{y\rightarrow +\infty }I_{f}(y)=\int\limits_{%
\mathbb{R}^{2g}}\,f(x)\,dx.
\]
\label{thm:mean_value}
\end{lem}

\textbf{proof}   Let $l=(m,n)^{T}\in \mathbb{Z}^{g}\times \mathbb{Z}^{g}$.
We now evaluate the integral $ I_{f}(y)$:
\[
 I_{f}(y)=\int\limits_{V_g}\text{ \ } \sum\limits_{\begin{subarray}{l}
  m,n \in \mathbb{Z}^g  \\
  (m,n) \ne (0,0)
\end{subarray}}  f\left(
\begin{array}{c}
y(m+Xn) \\
\frac{1}{y}n
\end{array}
\right) \text{ \ }d X \text{.}
\]

Here the sum and the integral are interchangeable. We split the sum over
$n$ into the part where $n=0$ and the part where $n\neq 0$. This way we
obtain for $I_{f}(y)$:
\[
 I_{f}(y)=\int\limits_{V_g}\text{ \ }\sum_{n\mathbb{=}0}\sum_{m\in \mathbb{Z}^{g}\backslash \{0\}}f\left(
\begin{array}{c}
y(m+Xn) \\
\frac{1}{y}n
\end{array}
\right) \text{ \ }d X+\sum_{n\in \mathbb{Z}^{g}\backslash \{0\}}c_{n}(X)  \text{.}
\]
For the first summand we obtain:
\[
\int\limits_{V_g}\ \ \sum_{n\mathbb{=}0}\sum_{m\in \mathbb{Z}^{g}\backslash
\{0\}}f\left(
\begin{array}{c}
y\cdot m \\
0
\end{array}
\right) \ \ d X= \sum_{m\in \mathbb{Z}^{g}\backslash
\{0\}}f\left(
\begin{array}{c}
y\cdot m \\
0
\end{array}
\right) \text{.}
\]

In the second summand we sum over $n\in \mathbb{Z}^{g}\backslash \{0\}$.
To structure this part of the proof, we first have to fix some notations
concerning the matrix $X$.\\
Let $M$ be the set of index-pairs
\[
M=\{(i',j') \mid (i'\leq j') \wedge (i'\leq g+1-j') \}.
\]
For $(i,j) \in M$ we have
\begin{eqnarray}
\label{equ:ij}
X_{i,j}=X_{j,i}&=&\sigma(i,j) X_{i,g+1-j}=\sigma(i,j) X_{g+1-j,i}, \text{ \ where}\\
\sigma (i,j) &=& \left\{ {\begin{array}{*{20}c}
   {-1}  \\
   { 1}  \\
 \end{array} } \right. \text{ \ if}\begin{array}{*{20}c}
   {i + j \text{ \ even}}  \\
   {i + j \text{ \ odd}}  \\
 \end{array}  \text{.} \nonumber 
\end{eqnarray}

For the proof we will have to consider pairs of rows and columns of the matrix $X$ of $P_Z$ simultaneously. We call $X_1,...,X_g$ the rows of the matrix $X$ and use the same notation for $P_Z$.
For the proof, we will have to consider the pairs $X_i$ and $X_{g+1-i}$, respectively. To abbreviate the notations we denote
\[
h(i)=g+1-i, \text{ \ for all \ } i \in \{1,...,g\}.
\]

Considering the terms of $n_{1}$ and $n_{h(1)}=n_g$ in the sum $\sum_{n\in \mathbb{Z}^{g}\backslash \{0\}}c_{n}(X)$ above, one has:
\begin{eqnarray*}
\sum_{n\in \mathbb{Z}^{g}\backslash \{0\}}c_{n}(X) &=&\sum_{(n_{1}^{\surd
},n_2,..,n_{g}^{\surd
})\backslash \{0\}}(\sum_{n_{1}=0 \wedge n_{g}=0}\int\limits_{V_g}%
\text{ \ }\sum_{m\in \mathbb{Z}^{g}}f\left(
\begin{array}{c}
y(m+Xn) \\
\frac{1}{y}n
\end{array}
\right) \text{ \ }d X)+ \\
&&\sum_{(n_{1}^{\surd },n_2...,n_{g}^{\surd })}(\sum_{n_{1}\neq
0 \vee  n_g  \neq 0} \int\limits_{[0,1]^{v_g-2}}\text{\ }\sum_{(m_1^{\surd },m_{2},...,m_{g}^{\surd })%
}F_{1,1}(m,n,X)\text{ \ }d X\backslash (X_{1,1},X_{1,g}) ) \text{.}
\end{eqnarray*}

Here $F_{1,1}(m,n,X)$ is defined as:
\[
F_{1,1}(m,n,X)=\sum_{m_{1},m_g}\int\limits_{[0,1]^2}\text{ \ }f\left(
\begin{array}{c}
y(m+Xn) \\
\frac{1}{y}n
\end{array}
\right) \text{ \ }d X_{1,1} d X_{1,g} \text{.}
\]
Analogously we define $\{F_{i,j}(m,n,X)\}_{(i,j) \in M}$ by
\[
F_{i,j}(m,n,X)=\sum_{m_{i},m_{h(i)}}\int\limits_{[0,1]^2}\text{ \ }f\left(
\begin{array}{c}
y(m+Xn) \\
\frac{1}{y}n
\end{array}
\right) \text{ \ }d X_{i,j}d X_{i,h(i)}\text{.}
\]

We now evaluate $F_{1,1}(m,n,X)$. We consider only the first and gth entry of $P_Z \cdot \left( {\begin{array}{*{20}c}
   {m}  \\
   {n}  \\
 \end{array} } \right)$ and note that $X_{1,1}$ and $X_{1,g}$ only occur in these two entries. $F_{1,1}(m,n,X)$ can now be written as
\[
F_{1,1}(m,n,X)=\sum_{m_{1},m_g}\int\limits_{[0,1]^2}\text{ \ }f\left(
\begin{array}{c}
y\cdot (m_{1}+\sum\limits_{j\neq 1,j \neq g}X_{1,j}n_{j}+X_{1,1}n_{1}+X_{1,g}n_g) \\
...\\
y\cdot (m_{g}+\sum\limits_{j\neq 1,j \neq g}X_{g,j}n_{j}-X_{1,1}n_{g}+X_{1,g}n_1)
\end{array}
\right) \text{ \ }d X_{1,1}d X_{1,g}\text{,}
\]
due to the relations in (\ref{equ:ij}).

These two rows can be written as
\begin{equation}
y\left( {A_{1g} \left( {\begin{array}{*{20}c}
   {X_{1,1} }  \\
   {X_{1,g} }  \\
 \end{array} } \right) + \left( {\begin{array}{*{20}c}
   {m_1 }  \\
   {m_g }  \\
 \end{array} } \right) + \left( {\begin{array}{*{20}c}
   {\lambda _{11} }  \\
   {\lambda _{1g} }  \\
 \end{array} } \right)} \right),
\label{eqn:A}
\end{equation}

where $A_{1g}$ is the matrix
$ \left( {\begin{array}{*{20}c}
   {n_1 } & {n_g }  \\
   {-n_g } & {n_1 }  \\
 \end{array} } \right)$. As $n_1 \neq 0$ or  $n_g\neq 0$ , $\det(A_{1g})= n_{1}^2 + n_{g}^2 \ne 0$ and we apply the  two dimensional formula for integration by substitution to these two rows. \\
 Note that $A_{1g} \cdot [0,1]^2$ is a parallelogram. Furthermore the translates of the parallelogram $A_{1g} \cdot [0,1]^2$ by $(m_1,m_g) \in \mathbb{Z}^2$ cover disjointly $\det(A_{1g})$ copies of $\mathbb{R}^2$. By integration by substitution we therefore obtain for $F_{1,1}(m,n,X)$:

\begin{equation*}
F_{1,1}(m,n,X) = \frac{\det(A_{1g})}{\det(A_{1g})} \int\limits_{\mathbb{R}^2}\text{ \ }f\left(
\begin{array}{c}
y\cdot X_{1,1} \\
... \\
y \cdot X_{1,g} \\
\end{array}
\right) \text{ \ }d X_{1,1} d X_{1,g} \text{.}
\end{equation*}

We proceed the same way with $X_{1,2}=X_{2,1},...,X_{1,\frac{g}{2}}=X_{\frac{g}{2},1}$ by successively integrating
over $F_{1i}(m,n,X)$. In each step the variables $X_{1,i}$ and $X_{1,h(i)}$ occur only twice and we use again integration by substitution. Here the determinants of the corresponding transformation matrices $A_{1i}$ are strictly positive, due to the shape of the matrix $X$.\\
Therefore we obtain by successive integration:

\begin{eqnarray}
&&\sum_{(n_{1}^{\surd },n_2...,n_{g}^{\surd })}\sum_{n_{1}\neq
0 \vee n_g \neq 0}\int\limits_{[0,1]^{v_{g}-2}}\text{ \ }\sum_{(m_{1}^{\surd},m_2...,m_{g}^{\surd})}F_{11}(m,n,X)\text{ \ }d X\backslash (X_{1,1},X_{1,g})  \nonumber \\
&=&\sum_{(n_{1}^{\surd },n_2...,n_{g}^{\surd })}(\sum_{n_{1}\neq
0 \vee n_g \neq 0}\int\limits_{[0,1]^{v_{g}-g}}(\frac{1}{y^{g}}\int\limits_{\mathbb{R}^{g}}\text{ \ }f\left(
\begin{array}{c}
X_{1} \\
\frac{1}{y}n
\end{array}
\right) \text{ \ }d X_{1})\text{ \ }d X\backslash X_{1})
\nonumber \\
&=&\frac{1}{y^{g}}\sum_{(n_{1}^{\surd },n_2...,n_{g}^{\surd })}(\sum_{n_{1}\neq 0 \vee n_g \neq 0}\int\limits_{\mathbb{R}^{g}}f\left(
\begin{array}{c}
X_{1} \\
\frac{1}{y}n
\end{array}
\right) \text{ \ }d X_{1})  \text{.}
\label{eq:succ_int}
\end{eqnarray}
Here $X_{1}=(X_{11},...,X_{1g})^{T}$. For the set $\mathbb{Z}^{g}\backslash \{0\}$ we obtain:
\[
\mathbb{Z}^{g}\backslash \{0\}= \\
\bigcup\limits_{i=1}^{g/2}%
\{n \in \mathbb{Z}^{g}\mid (n_{i}\neq 0)  \vee (n_{h(i)} \neq 0) \text{
and }n_{j}=0 \,\, \forall j \text {, such that } (j<i) \vee (j>h(i))\} \text{.}
\]

We now divide the set $\mathbb{Z}^{g}\backslash \{0\}$ successively into disjoint
unions of these subsets and sum up these subsets after integration over the $\{F_{i,j}(m,n,X)\}_{(i,j) \in M}$.
As in equation (\ref{eq:succ_int}) we obtain that

\[
\sum_{n\in \mathbb{Z}^{g}\backslash \{0\}}c_{n}(X)=%
\sum_{n\in \mathbb{Z}^{g}\backslash \{0\}}y^{-g}\int\limits_{\mathbb{R}^{g}}\text{
\ }f\left(
\begin{array}{c}
t \\
\frac{n}{y}
\end{array}
\right) \text{ \ }d t\text{.}
\]
As $f$ is continuous of compact support we obtain in total:
\begin{eqnarray*}
\lim\limits_{y\rightarrow +\infty }\text{ } I_{f}(y)
&=&\lim\limits_{y\rightarrow +\infty }(\sum_{n\in \mathbb{Z}^{g}\backslash
\{0\}}c_{n}(X)+ \sum_{m\in \mathbb{Z}^{g}\backslash \{0\}}f\left(
\begin{array}{c}
y\cdot m \\
0
\end{array}
\right) ) \\
&=&\text{ }\lim\limits_{y\rightarrow +\infty }(%
\sum_{n\in \mathbb{Z}^{g}\backslash \{0\}}y^{-g}\int\limits_{\mathbb{R}^{g}}\text{
\ }f\left(
\begin{array}{c}
t \\
\frac{n}{y}
\end{array}
\right) \text{ \ }d t) \\
&&+\text{ }\lim\limits_{y\rightarrow +\infty }\sum_{m\in \mathbb{Z}%
^{g}\backslash \{0\}}f\left(
\begin{array}{c}
y\cdot m \\
0
\end{array}
\right) \\
&=&\int\limits_{\mathbb{R}^{2g}}\text{ \ }f(x)\text{ \ }%
d x+0 = \int\limits_{\mathbb{R}^{2g}}\text{ \ }f(x)\text{ \ }%
d x\text{,}
\end{eqnarray*}
as for a continuous function $h:\mathbb{R}^{g}\rightarrow \mathbb{R}$ of compact support
one has
$$\lim\limits_{y\rightarrow+\infty }\sum_{n\in \mathbb{Z}^{g}\backslash \{0\}}y^{-g}\cdot h\left(\frac{n}{y}%
\right)=\int\limits_{\mathbb{R}^{g}} h\left( t\right) d t,$$
by integration over Riemann sums.   $\square$\\
\\
\section{Families of highly symmetric symplectic lattices}

We first present a construction that enables us to produce families of lattices with a large automorphism group
in any dimension $g$. These have always at least $g-1$ symmetries.\\
Consider the matrix $C_g \in M_g(\R)$: 
\begin{equation*}
C_g := 
\left(
\begin{array}{cc}
0 & Id_{g-1} \\
1 & 0
\end{array}
\right) \text{, \ \ where }C_g^{T}=C_g^{-1}\text{.}
\end{equation*}
It is easy to see that $(C_g)^g=Id_g$ and the multiplicative matrix group $\mathcal{Z}_g$ generated by $C_g$ is isomorphic to $(\Z_g,+)$: 
\[
(\mathcal{Z}_g,\cdot) \simeq (\Z_g,+).
\]
The group $(\mathcal{Z}_g,\cdot) \simeq (\Z_g,+)$ acts by conjugation on $M_{g}(\mathbb{R)}$. We denote by $\sta(\Zc)$ its fixed points in $M_{g}(\mathbb{R)}$:
\[
\sta(\Zc)=\{A\in M_{g}(\mathbb{R)} \mid (C_g)^{-1}AC_g=A\}.
\]
It follows from this fixpoint equation that $\sta(\Zc)$ consists exactly of the circulant matrices (see \cite{ab}, \textbf{Theorem 2.1} for a proof). A characterization of the circulant matrices and the example further below is given in \cite{ch} and \cite{di} in terms of an eigenvalue decomposition.\\
Now consider a matrix $A \in \sta(\Zc) \cap GL_g(\R)$ and the corresponding lattice $A\cdot \Z^g$. As $C_g $ is in  $O_g(\F_2)$ and therefore in $O_g(\R) \cap GL_g(\Z)$, we can set $O=(C_g)^{-1}$ and $R=C_g$ in  equation (\ref{eqn:sym}). It follows that the lattice $A\cdot \Z^g$ has $g-1$ symmetries. In other words the lattices, whose basis can be represented by a circulant matrix have at least $g-1$ disjoint symmetries. \\
Let $\otimes$ denote the tensor product. As $\mathcal{Z}_g$ is a multiplicative group, we have 
\[
     \mathcal{Z}_{g_1} \otimes  \mathcal{Z}_{g_2}  \simeq   \Z_{g_1} \oplus  \Z_{g_2}.  
\]
By the fundamental theorem of finite abelian groups we can therefore represent any finite abelian group
of order $g$ as a multiplicative matrix group, which is a subgroup of $O_g(\R)$. We can find representations of these groups using the Kronecker product on the generator matrices.  \\
\\
In the following we will examine the most simple of these examples, 
\[
\mathcal{Z}^n_{2}= \mathop{\otimes} \limits_{i=1}^{n} \mathcal{Z}_2.
\]
We will first characterize the matrices that constitute the fixed points of this group under conjugation. These have interesting algebraic and geometric properties. We will then shortly discuss the connection with the Barnes-Wall lattices. Then  we will establish a family of highly symmetric symplectic lattices with the help of this matrix group.\\
To understand the geometry of these matrices consider the following construction.
Let $J_g \in O_g(\mathbb{R})$ be the matrix

\begin{equation*}
J_g=\left(
\begin{array}{cccc}
0 & ... & 0 & 1 \\
... &  & 1 & 0 \\
0 & 1 & ... &   \\
1 & 0 & ... & 0 \\
\end{array}
\right) \text{, \ where }J_g=J_g^{T}=J_g^{-1}\text{.}
\end{equation*}

For $A \in M_g(\R)$, we denote by $A_{T}$ the matrix obtained by reflecting $A$ along its second main diagonal. We have
\begin{equation}
J_g A J_g = A\Leftrightarrow A=A^T =A_T  \text{.}
\label{eq:Jsym}
\end{equation}

We have $J_2 = C_2$ and $(\{Id_2,J_2\},\cdot) = (\mathcal{Z}_{2},\cdot)$. Using the Kronecker product $\otimes$ on $J_2$, we can find a set of generating matrices for $\mathcal{Z}^n_{2}$. For $n \in \N$ and $g=2^n$ set 
\[
   J_{2^n}^{n}= J_2^{\otimes n} = J_{2^n}  ,J_{2^n}^{n-1}= J_2^{\otimes n-1} \otimes Id_2
   ,...,J_{2^n}^{1}=J_2 \otimes Id_2^{\otimes n-1}.
\]
We have $J_{2^n}^{k}={J_{2^n}^{k}}^{-1}={J_{2^n}^{k}}^{T}$. These matrices generate $\mathcal{Z}^n_{2}$. Denote by
\[
\sta(\mathcal{Z}^n_2)=\{A \in M_{2^n}(\mathbb{R)} \mid J_{2^n}^{k}AJ_{2^n}^{k}=A , \text{ \ \ for \ \ } k \in \{1,...,n\}\}
\]
the set of fixed points in $M_{2^n}(\mathbb{R})$ under the action of $\mathcal{Z}^n_2$. \\
\\
As in the case of $J_{2^n}=J_{2^n}^{n}$ (see (\ref{eq:Jsym})) the conjugation with the other generating matrices $(J_{2^n}^{k})_{k=1,...,n-1}$ of $\mathcal{Z}^n_2$ imply further symmetries for the matrices in $\sta(\mathcal{Z}^n_2)$. 

\begin{figure}[h!]
\SetLabels
\L(.28*.10) $g=2$\\
\L(.47*.10) $g=4$\\
\L(.66*.10) $g=8$\\
\endSetLabels
\AffixLabels{%
\centerline{%
\includegraphics[height=5cm,width=10cm]{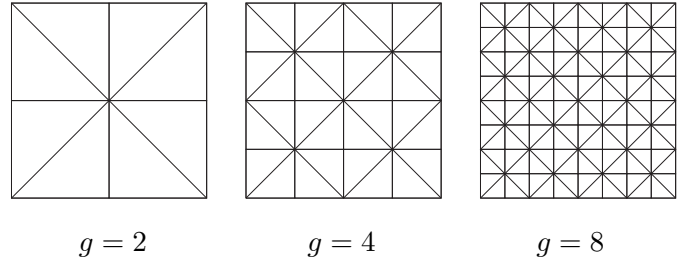}}}
\caption{Symmetries of the matrices in $\sta(\mathcal{Z}^n_2) \subset M_g(\R) $ for $g=2,4$ and $8$.}
\label{fig:sym_A}
\end{figure}

The symmetries among the matrix entries of $\sta(\mathcal{Z}^n_2) \subset M_{g}(\R) $ for $g=2,4$ and $8$ are depicted in Figure \ref{fig:sym_A}. Here each box represents a matrix entry and each diagonal line an axis of symmetry. These symmetries equally define $\sta(\mathcal{Z}^n_2)$. \\

We obtain another inductive description of $\sta(\mathcal{Z}^n_2)$ in the following way. Let $\mathbb{A}_2$ be the set of matrices
\[
\mathbb{A}_2=\{C \in M_{2}(\mathbb{R}) \mid C = \left( {\begin{array}{*{20}c}
   {a } & {b }  \\
   {b } & {a }  \\
 \end{array} } \right), \text{ \ \ } a,b \in \mathbb{R}\}.
\]
We define inductively $\mathbb{A}_{2^{n+1}}$ to be the set of matrices
\begin{equation}
\mathbb{A}_{2^{n+1}}=\{C \in M_{2^{n+1}}(\mathbb{R}) \mid 
C = \left( {\begin{array}{*{20}c}
   { A} & {B}   \\
   { B}   & {A}  \\
 \end{array} }  \right), \text{ \ such that } A,B  \in \mathbb{A}_{2^{n}}\}.
\label{eq:A_ind}
\end{equation}
A matrix in $\mathbb{A}_{2^{n}}$  is entirely determined by its first row. It is easy to prove by induction that
\begin{equation*}
A \in \mathbb{A}_{2^n} \Leftrightarrow J_{2^n}^{k} A J_{2^n}^{k}=A \text{\ \ for all } k \in \{1,..,n\} \Leftrightarrow A \in \sta(\mathcal{Z}^n_2).
\label{eq:comm}\\
\end{equation*}
Hence $\mathbb{A}_{2^n} = \sta(\mathcal{Z}^n_2)$. For these matrices we have the following eigenvalue decompositon:   

\begin{lem}
Let $C \in M_g(\R)$ be a matrix in $\sta(\mathcal{Z}^n_2)$. Then the columns of the matrix 
\[
V_g=\left( {\begin{array}{*{20}c}
   {1 } & {1 }  \\
   {1 } & {-1 }  \\
 \end{array} } \right)^{\otimes n}    , \text{ \ \ where \ \ } (V_g)^{-1} = \frac{1}{g} \cdot V_g  
\] 
form an orthogonal basis of eigenvectors. If $v_i$ is an eigenvector and $c^1=(c_{11},...,c_{1g})$ the first row of a matrix $C \in \sta(\mathcal{Z}^n_2)$ then $d_i = c^1 \cdot v_i$ is the corresponding eigenvalue for $v_i$.\\
\end{lem}
This decomposition shows that these matrices belong to the patterned matrices described in \cite{ch} or \cite{di}.\\

\textbf{proof}  The proof follows by induction: \\
We use the characterization of $\sta(\mathcal{Z}^n_2)= \mathbb{A}_{2^n}$ in equation (\ref{eq:A_ind}). A simple calculation shows that the statement is true for the matrices in $\mathbb{A}_{2}$. We now assume that the eigenvectors and eigenvalues of $\mathbb{A}_{2^{n}}$ have the form stated above. Then for any matrix  $C \in \mathbb{A}_{2^{n+1}}$, we have for some $A,B \in \mathbb{A}_{2^{n}}$ by our induction hypothesis:  
\[  
    C = \left({\begin{array}{*{20}c}
   { A} & {B}   \\
   { B}   & {A}  \\
 \end{array} }  \right) = 
   \left({\begin{array}{*{20}c}
   { V_{2^{n}}} & {0 }  \\
   { 0}   & {V_{2^{n}}}  \\
 \end{array} }  \right) 
  \left({\begin{array}{*{20}c}
 { D_A } & {D_B }  \\
   { D_B}   & {D_A}  \\
 \end{array} }  \right) 
  \left({\begin{array}{*{20}c}
   { V_{2^{n}}} & {0 }  \\
   { 0}   & {V_{2^{n}}}  \\
 \end{array} }  \right)^{-1}. 
\]
Here $D_A$ and $D_B$ are the diagonal matrices that contain the eigenvalues of $A$ and $B$, respectively. Defining $D_C$ in the same manner with respect to $C$, we have: 
\[ 
\left({\begin{array}{*{20}c}
 { Id_{2^n} } & {Id_{2^n} }  \\
   { Id_{2^n} }   & {-Id_{2^n}}  \\
 \end{array} }  \right)
    \left({\begin{array}{*{20}c}
 { D_A } & {D_B }  \\
   { D_B}   & {D_A}  \\
 \end{array} }  \right) 
\left({\begin{array}{*{20}c} 
  { Id_{2^n} } & {Id_{2^n} }  \\
   { Id_{2^n}}   & {-Id_{2^n}}  \\
 \end{array} }  \right)^{-1}
  = D_C.
\]
Combining these two equations we obtain our statement by induction. This completes our proof.  \ \ \     $\square$ 
\\  
  
Now consider the matrices in $\sta(\mathcal{Z}^n_2) \cap GL_g(\R)$. The columns of such a matrix form a basis of a $g$-dimensional lattice. As in the case of the invertible circulant matrices it follows that the corresponding lattices, whose basis can be represented by such a matrix have at least $g-1$ disjoint symmetries. \\
\\
Summarizing these results we have shown:
\begin{thm}
The matrix group $(\mathcal{Z}^n_2,\cdot) \simeq ({\Z_2}^n,+)$ acts by conjugation on $M_{2^n}(\R)$. Let $\sta(\mathcal{Z}^n_2)$ be the set of fixed points under the group action. Then $\sta(\mathcal{Z}^n_2)=\mathbb{A}_{2^{n}}$ and the set of fixed points can be characterized by its eigenvalue decomposition. Every lattice $\Lambda = A\cdot \mathbb{Z}^{2^n}$, where $A \in \sta(\mathcal{Z}^n_2) \cap GL_{{2^n}}(\R)$ has at least $2^{n}-1$ disjoint symmetries.
\end{thm}

An important class of lattices that occur exactly in dimensions of power of two are the Barnes-Wall lattices (see \cite{bw}). We could not prove or disprove that these lattices have a basis that has a matrix representation in $\sta(\mathcal{Z}^n_2)$. However we can always find a basis representation that shares some of the matrix symmetries with this set. To this end we use the definition of the Barnes-Wall lattices $BW_{2^n}$ as $2^n$ dimensional
lattices over the Gaussian integers $\Z + i\Z$:
\begin{defi}
For $n>1$, $BW^{\C}_{2^n}$ is the $2^n$ dimensional lattice over $\Z + i\Z$ generated
by the columns of the n-fold Kronecker product
\[
BW^{\C}_{2^n}=\left( {\begin{array}{*{20}c}
           {1 } & {1 }  \\
           {1 } & {i }  \\
           \end{array} } \right)^{\otimes n}.
\] 
\end{defi}

\begin{figure}[h!]
\SetLabels
\L(.28*.10) $g=4$\\
\L(.47*.10) $g=8$\\
\L(.66*.10) $g=16$\\
\endSetLabels
\AffixLabels{%
\centerline{%
\includegraphics[height=5cm,width=10cm]{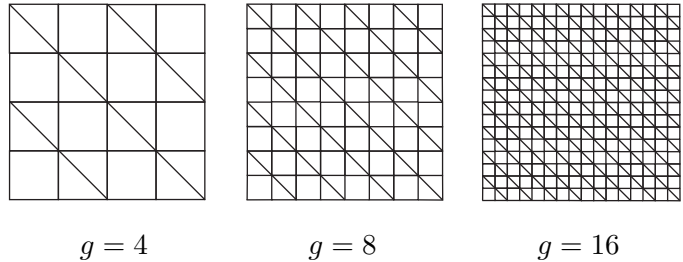}}}
\caption{Symmetries of the matrices in $BW'_g$ for $g=4,8$ and $16$.}
\label{fig:sym_BW}
\end{figure}

This definition is based on \cite{ni}. In \cite{ni} another matrix of basis vectors is used. The matrix of our definition can be obtained by a simple change of basis.  Using the identification 
\[
    a+ib \mapsto    
\left( {\begin{array}{*{20}c}
   {a } & {b }  \\
   {-b } & {a }  \\
 \end{array} } \right)   
\] 
we obtain a matrix $BW^{\R}_{2^{n+1}}$ in $\R^{2^{n+1}}$, whose columns form a basis of minimal vectors for the Barnes-Wall lattices. Now 

\[
\left( {\begin{array}{*{20}c}
   {a} & {b }  \\
   {-b } & {a }  \\
 \end{array} } \right) \cdot  J_2 =  
 \left( {\begin{array}{*{20}c}
   {b } & {a }  \\
   {a } & {-b }  \\
 \end{array} } \right).
\] 
Applying a change of basis, we obtain another basis $BW'_{2^{n+1}} =   BW^{\R}_{2^{n+1}} \cdot J_2 \otimes Id_2^{\otimes n}$. These matrices $BW'_{2^{n+1}}$  have the symmetries depicted in Figure \ref{fig:sym_BW}. 
Here each box represents a matrix entry and each diagonal line an axis of symmetry.\\

After scaling these lattices to determinant $1$, we obtain that the squared norm of the minimal vectors is $\sqrt{\frac{g}{2}}$. The scaled versions of the Barnes-Wall lattices are alternatively 2-modular and unimodular and the unimodular lattices are symplectic (see \cite{ber1}, section 7).\\
\\
We now apply the construction of symmetric matrices and lattices to the moduli space of symplectic lattices of dimension $2^n$.
Remember that for  $Z=X+i\cdot Y \in \mathcal{H}_{g}$ we have

\begin{equation*}
P_{Z}=P_{X,Y}=\left(
\begin{array}{cc}
\sqrt{Y^{-1}} & 0 \\
0 & \sqrt{Y}
\end{array}
\right) \left(
\begin{array}{cc}
Id_g & X \\
0 & Id_g
\end{array}
\right) \in M_{2g}(\mathbb{R)}\text{.}
\label{eqn:PZH2}
\end{equation*}

We now consider the following set of lattices defined by a $Z=X+i\cdot Y \in \mathcal{H}_{2^n}$:
\begin{equation}
\mathbb{A}^{s}_{2^n}= \{Z=X+i\cdot Y \in \mathcal{H}_{2^n} \mid \sqrt{Y},X \in \sta(\mathcal{Z}^n_2) \}.
\label{eqn:bsymp}
\end{equation}

For the elements of $\mathbb{A}^{s}_{2^n}$ we have:

\begin{thm} Let $g=2^n$. If $Z \in \mathbb{A}^{s}_{2^n}$, then the symplectic lattice $P_Z\cdot \Z^{2g}$ has at least $g-1$ symmetries.
\label{thm:K_g}
\end{thm}

\textbf{proof} To prove that $P_Z$ has $g-1$ symmetries, we note that for the generators $\left(J_{2^n}^{k}\right)_{k=1,..n}$ of $\mathcal{Z}^n_2$, we have for all $k \in \{1,..,n\}$:
\[
\left( {\begin{array}{*{20}c}
   J_{g}^{k} & 0  \\
   0 & J_{g}^{k}  \\
 \end{array} } \right)
 \left(
\begin{array}{cc}
\sqrt{Y^{-1}} & 0 \\
 0 & \sqrt{Y}\\
\end{array}
\right)
\left( {\begin{array}{*{20}c}
   J_{g}^{k} & 0  \\
   0 & J_{g}^{k}  \\
 \end{array} } \right)
=\left(
\begin{array}{cc}
\sqrt{Y^{-1}} & 0 \\
0 & \sqrt{Y} \\
\end{array}
\right) .
\]
Here we use that for $A \in \sta(\mathcal{Z}^n_2), A \in GL_{2^n}(\R) \Rightarrow A^{-1} \in \sta(\mathcal{Z}^n_2)$. Furthermore
\[
\left( {\begin{array}{*{20}c}
   J_{g}^{k} & 0  \\
   0 & J_{g}^{k}  \\
 \end{array} } \right)
 \left(
\begin{array}{cc}
Id_{g} & X \\
0 & Id_{g} \\
\end{array}
\right)
\left( {\begin{array}{*{20}c}
   J_{g}^{k} & 0  \\
   0 & J_{g}^{k}  \\
 \end{array} } \right)
=\left(
\begin{array}{cc}
Id_{g} & X \\
0 & Id_{g} \\
\end{array}
\right)  , \text{ \ \ for all } k \in \{1,..,n\}.
\]
It follows that
\[
\left( {\begin{array}{*{20}c}
   J_{g}^{k} & 0  \\
   0 & J_{g}^{k}  \\
 \end{array} } \right)
 \cdot
 P_Z
 \cdot
\left( {\begin{array}{*{20}c}
   J_{g}^{k} & 0  \\
   0 & J_{g}^{k}  \\
 \end{array} } \right)
=P_Z  \\
, \text{ \ \ for all } k \in \{1,..,n\},
\]
as $ J_{g}^{k} \cdot J_{g}^{k} = Id_{g} , \text{ \ \ for all } k \in \{1,..,n\}.$ \\

As this is true for the generators $(J_{g}^{k})_{k=1,..,n}$ of $\mathcal{Z}^n_2$, these equations are true for all elements of $\mathcal{Z}^n_2$.
Therefore $P_Z$ has the demanded symmetries. (see equation (\ref{eqn:sym})).  $\square$

\begin{con}
If $g=2^n$ and $\delta_{2g}$ Hermite's invariant for symplectic lattices of dimension $g$, then
\[
 \frac{1}{\pi }\sqrt[g]{2 g \cdot g!} \leq \delta_{2g}.
\]
\label{thm:con_her}
\end{con}
The idea is to use in \textbf{Theorem \ref{thm:Hermite_mu}} a subset of the symmetric lattices, whose matrix representation of a basis is in $\mathbb{A}^{s}_{2^n}$ (see (\ref{eqn:bsymp})). Due to the symmetries of these lattices the number of vectors of a certain length in a lattice $P_{X,\frac{1}{y^2}Id_{g}} \cdot \Z^g$ is always a multiple of $2g$. Therefore we could in principle apply the methods of the proof of \textbf{Theorem \ref{thm:Hermite_mu}} to prove the above conjecture. The problem is to evaluate the integral in the case, where a linear transformation as in equation (\ref{eqn:A}) is not possible. Here the matrix in the transformation is singular in some cases and we can not proceed.\\ 
Whereas this bound given in the conjecture is of order $g$, the symplectic Barnes-Wall lattices in dimension $g$ have minimal vectors whose squared length is $\sqrt{\frac{g}{2}}$.\\
\\
As already mentioned in the introduction, a better lower bound for Hermite's constant $\gamma_g$ than in \textbf{Theorem \ref{thm:Hermite_mu}} was given by Rogers in \cite{ro} for real lattices of dimension $g$.

\section*{Acknowledgments}
The author has been supported by the Swiss National Science Foundation (grant nr. 200021 - 125324). I would like to thank Peter Buser for his support and encouragement. I would also like to thank Greg Kuperberg for pointing me to the given construction of the Barnes-Wall lattices on \textit{mathoverflow}.

\end{document}